\documentclass[36pt,reqno]{amsart}
 \usepackage[bookmarks]{hyperref}

\newcommand\LL{\mathcal L}
\newcommand\OO{\mathcal O}
\newcommand\E{\mathcal E}
\newcommand\PP{\mathbb P}

\newcommand\X{\mathcal X}
\newcommand\Y{\mathcal Y}

\newcommand{\Hilb}{\operatorname{Hilb}}
\newcommand{\LCH}{\lambda_{\operatorname{CH}}}
\newcommand{\LCM}{\lambda_{\operatorname{CM}}}
\newcommand{\LHilb}{\lambda_{\operatorname{Hilb}}}

\makeatletter \@addtoreset{equation}{section} \makeatother

\newtheorem{thm}[equation]{Theorem}
\newtheorem{lem}[equation]{Lemma}
\newtheorem{cor}[equation]{Corollary}
\newtheorem{prop}[equation]{Proposition}

\theoremstyle{definition}
\newtheorem{defn}[equation]{Definition}

\newtheorem{example}[equation]{Example}
\newtheorem{rmk}[equation]{Remark}

\title[A note on Positivity of the CM line]{A note on Positivity of the CM line bundle}
\author{J. Fine and J. Ross}

\begin{document}
\bibliographystyle{abbrv}

\begin{abstract} We show that positivity of the CM line associated to a family of polarised varieties is intimately related to the stability of its members. We prove that the CM line is nef on any curve which meets the stable locus, and that it is pseudoeffective (i.e.\ in the closure of the effective cone) as long as there is at least one stable fibre.  We give examples showing that the CM line can be strictly negative or strictly positive on curves in the unstable locus. 
% We also compare
%  this to the work of Fujiki-Shumacher.
\end{abstract}

\maketitle

\section{Introduction}

A famous conjecture of Yau relates the existence of K\"ahler metrics
of constant scalar curvature to stability in Geometric Invariant
Theory (GIT).  It is thought that this stability notion should be with
respect to the CM line bundle on the Hilbert scheme, originally
defined by Tian \cite{tian(94):k_energ_hyper_stabil}.  Unfortunately,
it is not clear directly from the definition when this line bundle is positive (e.g.\ ample), so one cannot define stability in terms of the
existence of non-vanishing invariant sections as is standard in GIT.
Instead, the definition of stability that is commonplace is made by analogy with the Hilbert-Mumford criterion, and requires that for all one parameter subgroups the the Hilbert-Mumford weight function has the favourable sign.

We start this note by defining the CM line following Paul--Tian
\cite{paul_tian:algeb_analy_k_stabiI} and then show how it arises as the
leading order term in an expansion of some naturally defined line
bundles.  We then turn to positivity properties of the CM line (we use
the term positivity loosely here to mean ample, nef, effective and so
on). Roughly speaking, we show that the CM line is non-negative on the ``stable'' locus (more precisely, the asymptotically Hilbert semistable locus) and that, as long as the stable locus is non-empty,
the CM line is pseudoeffective.  We end with examples showing that outside the stable locus the CM line is badly behaved: there are examples of families of smooth unstable surfaces for which the CM line is strictly negative, and other such examples for which it is positive.
\vspace{4mm}

{\bf Notation and Conventions: } If $\LL$ is a line bundle its powers
are $\LL^k=\LL^{\otimes k}$ and $\LL^{-k} = (\LL^*)^{\otimes k}$.  If $\X\to
B$ is a family of schemes the fibre over a point $b\in B$ is $\X_b$
and if $\LL\to \X$ is a line bundle then $\LL_b = \LL|_{\X_b}$.  By a
polarised variety or scheme $(X,L)$ we mean a choice of ample line
bundle $L$ on $X$.

We say that a line bundle on $B$ is \emph{nef} if it has non-negative
degree on any irreducible curve $C$ in $B$.  It is effective if some
positive tensor power has a section, and it is \emph{pseudoeffective}
if its first Chern class lies on the closure of the cone of effective
divisors in the N\'eron-Severi space $N^1(B)$.  We work throughout
over the complex numbers.  \vspace{4mm}

{\bf Acknowledgements: } We would like to thank Sean Paul and Gang
Tian for useful conversations, and explaining to us their definition of the CM line \cite{paul_tian:algeb_analy_k_stabiI}.

\section{Definition of the CM line bundle}

Let $\pi\colon \X\to B$ be a proper flat morphism of schemes of constant relative dimension $n\ge 1$ and let $\LL$ be a relatively ample line bundle on $\X$.  We will assume throughout that
$B$ is irreducible, and that $\X$ has pure dimension.  

The CM line is defined in terms of the determinant of the pushdown of $\LL^k$ given by
$$\det \pi_!(\LL^k) = \det R^{ \bullet}\pi(\LL^k) = \bigotimes_{i} \left( \det R^i\pi_*(\LL^k)\right)^{(-1)^i}.$$
As $\LL$ is relatively ample the terms $R^i\pi_*(\LL^k)$ vanish for
$i>0$ and $k\gg 0$ so $\det \pi_!(\LL^k) = \det \pi_*(\LL^k)$.  We will
rely on the fact that $\det \pi_!(\LL^k)$ has a polynomial expansion in
terms of some fixed line bundles $\lambda_i$ on $B$.

\begin{thm}[Mumford--Knudsen]\label{thm:mumford}There exist line bundles $\lambda_{i}=\lambda_i(\X,\LL)$ on $B$
  such that for all $k$ 
\begin{equation}\label{mumfordexpansion}
\det \pi_!(\LL^k) \cong \lambda_{n+1}^{\binom{k}{n+1}} \otimes \lambda_{n}^{\binom{k}{n}} \otimes \cdots \otimes \lambda_0.
\end{equation}
%Moreover the $\lambda_i$ as well as the isomorphism are canonical and commute with base change.
\end{thm}

\begin{proof}
  See (GIT \cite{mum94:geometric} page 230) or 
  (Theorem 4
  \cite{knudsen_mumford(76):projec_modul_space_stabl_curves}).
%In that
%  theorem pick $r\gg 0$ so $\LL^r$ is relatively very ample and if
%  $E^*=\pi_*(\LL^r)$ then $E$ is locally free and $\X\subset \PP(\pi_* E^*)$.
%  The complex $\mathcal{F}^{\bullet}$ is simply $\OO_{\X}$.
%  which satisfies the technical condition $Q_{(n)}$
%  (\cite{knudsen_mumford(76):projec_modul_space_stabl_curves} p.50) by
%  flatness of $\pi$.
\end{proof}

We now come to the definition of the CM line.  Although it is not
immediately apparent, it turns out that this is a very natural line
bundle to consider (see Section \ref{sec:expansion}).  Write the
Hilbert polynomial of the fibres of $\X$ as
$$p(k) = \chi(\LL_b^k)= a_0k^n + a_1k^{n-1} + O(k^{n-2}),$$ 
which by flatness of the family is independent of $b\in B$.  As $\LL_b$
is ample the term $a_0$ is strictly positive, so we
can set
$$\mu = \mu(\X,\LL)=\frac{2a_1}{a_0}.$$

\begin{defn}(Paul--Tian
  \cite{paul_tian:algeb_analy_k_stabiI}) \label{def:cmline}The
  \emph{CM line bundle} associated to the family $(\X,\LL)$ is
  $$
  \LCM = \LCM(\X,\LL)=\lambda_{n+1}^{\mu +n(n+1)} \otimes
  \lambda_{n}^{-2(n+1)}.$$
\end{defn} \vspace {2mm}

%\begin{rmk}
%  It is important to note that the $\lambda_i$ are not canonical, as
%  they depend on the choice of $r$ to make $\L^r$ relatively very
%  ample (i.e\ they depend on the embedding of $\X$ in $\PP(\pi_* L)$.
%  We will show however in Section 2 that the CM line is independent
%  of this choice. **is this correct?**%
%\end{rmk}

Paul--Tian show in \cite{paul_tian:algeb_analy_k_stabiI} that $\LCM$
agrees with Tian's original definition of the CM line in terms of the
pushdown of a certain virtual bundle on $\X$
\cite{tian(94):k_energ_hyper_stabil}.  For example, suppose $\pi\colon
\X\to B$ has a relative canonical bundle $K_{\X/B}$.  An easy
calculation with the Grothendieck--Riemann--Roch formula applied to
$\LL^k$ gives
    \begin{eqnarray*}
      c_1(\lambda_{n+1}) &=& \pi_* c_1(\LL)^{n+1}\\
      nc_1(\lambda_{n+1}) - 2c_1(\lambda_n) &=& \pi_* \left(c_1(\LL)^n c_1(K_{\X/B})\right).
    \end{eqnarray*}
Hence
\begin{equation}
c_1(\LCM) = \pi_* \left[\mu c_1(\LL)^{n+1} +(n+1)c_1(K_{\X/B})c_1(\LL)^n\right].\label{eq:cmchernclass}
\end{equation} \vspace {4mm}

The definition of the CM line is made so that the Hilbert--Mumford
weight function (in the sense of geometric invariant theory)
equals Donaldson's version of the Futaki invariant
\cite{donaldson(02):scalar_curvat_stabil_toric_variet} as in the
following lemma.  This is the reason for the link between $\LCM$ and
the problem of finding K\"ahler metrics of constant scalar curvature.
We refer the interested reader to \cite{tian00:canonical_kaehler} for
an introduction to this topic.

\begin{lem}[Paul--Tian]
  Suppose that $\mathbb C^\times$ acts on $(\X,\LL)$ covering an action
  on $B$, and that $b\in B$ is a fixed point.  Write the weight
  (i.e.\ the sum of the eigenvalues) of the induced action on
  $H^0(\X_b,\LL_b^k)$ as a polynomial in $k$ with coefficients
$$w(k) = b_0 k^{n+1} + b_1k^n + O(k^{n-1}).$$
Then the weight of the induced action on $\LCM|_b$ is the Futaki invariant
$$ F_1 = \frac{2(n+1)!}{a_0} \left( b_1a_0 -b_0a_1\right).$$
\end{lem}

\begin{rmk}\
 The CM line is homogeneous of order $n$, \emph{i.e.} if $r\in \mathbb N$ then 
$$\LCM(\X,\LL^r) = r^n\LCM(\X,\LL).$$
Thus by clearing denominators we can define $\LCM(\X,\LL)$ as a
$\mathbb Q$-line bundle when $\LL=\OO(D)$ for some relatively ample
$\mathbb Q$-divisor $D$, and it makes sense to talk about $\LCM$ being
positive (e.g.\ ample, nef or pseudoeffective) in this case.
\end{rmk}

\section{The CM line as a leading order term}\label{sec:expansion}

We now give an alternative description of the CM line in terms of some
naturally defined line bundles.  We retain the notation from the
previous section, so $\pi\colon \X\to B$ is a proper flat family of
schemes of relative dimension $n$, $\LL\to \X$ is a relatively ample
line bundle on $\X$ and the fibres of $\X$ have Hilbert
polynomial $p$.  We set
$$\lambda(k) = \det \pi_!(\LL^k).$$

Suppose that $\sigma$ is a line bundle on $B$.  If for $k\gg 0$ we
replace $\LL$ by $\LL\otimes \pi^* \sigma$ then by the projection
formula $\lambda(k)$ becomes
$$ \det \pi_!\left((\LL\otimes \pi^* \sigma)^k\right) = \det \left(\pi_*(\LL^k) \otimes \sigma^k\right) = \det \pi_*(\LL^k)\otimes \sigma^{kp(k)}.$$
Thus for $k\gg 0$ the line bundle
$$\LHilb(\X,\LL,k) = \lambda(k)^{p(1)} \otimes \lambda(1)^{-kp(k)}$$
is rigid, i.e.\ independent of twisting $\LL$ by $\pi^*\sigma$, and thus should be thought of as a natural line bundle associated purely to the family of polarised varieties $(\X_b,\LL_b)$.

To relate $\LHilb$ to $\LCM$ we replace $\LL$ by $\LL^r$ and consider
$$\LHilb(\X,\LL^r,k)= \lambda(kr)^{p(r)} \otimes \lambda(r)^{-kp(kr)}$$
which now depends on two variables $k$ and $r$.  The CM line is the
leading order term of this line bundle for $k\gg r\gg 0$ as in the
following proposition.

\begin{prop}\label{prop:leadingorder}
  Let $\X\to B$ and $\LL$ be as above.  Then
$$\LHilb(\X,\LL^r,k) = \LCM(\X,\LL)^{\frac{a_0}{2(n+1)!}k^{n+1} r^{2n}} \otimes \epsilon(r,k)$$
where $\epsilon(r,k)$ consists of ``lower
order terms'' in $k$ and $r$.  More precisely there exist fixed line
bundles $\epsilon_{i,j}$ on $B$ such that
$$\epsilon(r,k) = \bigotimes_{i,j} \epsilon_{i,j}^{k^i r^j}$$
where the product is over $i=0,\ldots,n+1$ and $j=0,\ldots, 2n+1$ with $j\le 2n-1$ if $i=n+1$.
\end{prop}
 
\begin{proof}
Recall from \eqref{thm:mumford} that there are line bundles $\lambda_i$ on $B$ and an expansion
\begin{equation}
\det \pi_!(\LL^k) \cong \lambda_{n+1}^{\binom{k}{n+1}} \otimes \lambda_{n}^{\binom{k}{n}} \otimes \cdots \otimes \lambda_0.
\end{equation}
For convenience, define the Cornalba--Harris line associated to
$(\X,\LL)$ as
   \begin{equation}
     \LCH=\LCH(\X,\LL)=   \lambda_{n+1}^{p(1)} \otimes \lambda(1)^{-a_0(n+1)!}.\label{eq:linecornalba}
   \end{equation}
This is the leading order piece of $\LHilb(\X,\LL,k)$ as $k\gg 0$ since, omitting terms of order $k^n$, we have
\begin{eqnarray*}
    \LHilb(\X,\LL,k) &=& \lambda(k)^{p(1)} \otimes \lambda(1)^{-kp(k)} \\
&=& \lambda_{n+1}^{\binom{k}{n+1}p(1)} \otimes \cdots \otimes \lambda(1)^{-k(a_0k^n + \cdots)} \\
&=& \left(\lambda_{n+1}^{k^{n+1} p(1)}\otimes \lambda(1)^{-a_0(n+1)! k^{n+1}} \right)^{\frac{1}{(n+1)!}} 
\otimes \cdots \\
&=& \LCH(\X,\LL)^{\frac{k^{n+1}}{(n+1)!} } \otimes \cdots.
\end{eqnarray*}
Replacing $\LL$ by $\LL^r$ this yields
\begin{equation}\label{eq:chexpansion}
  \LHilb(\X,\LL^r,k) =\LCH(\X,\LL^r) ^{\frac{k^{n+1}}{(n+1)!}} \otimes \cdots.
\end{equation}
 Thus we need to expand
$\LCH(\X,\LL^r)$ as a function of $r$.  To this end notice that
replacing $\LL$ by $\LL^r$ has the effect of replacing $\lambda_{n+1}$
by $\lambda_{n+1}^{r^{n+1}}$ and $a_0$ by $a_0r^n$.  So for $r\ge 1$
we have
\begin{eqnarray*}
\LCH(\X,\LL^r)&=&\lambda_{n+1}^{p(r)r^{n+1}} \otimes \lambda(r)^{-r^n a_0(n+1)!}.
\end{eqnarray*}
Now $\binom{r}{n+1} = \frac{ r^{n+1}}{(n+1)!} - \frac{n(n+1)}{2(n+1)!} r^n + O(r^{n-1})$, so  up to terms of order $r^{n-1}$,
$$\lambda(r) = \lambda_{n+1}^{ \binom{r}{n+1}} \otimes \lambda_n ^{\binom{r}{n}} \otimes \cdots = \lambda_{n+1}^{\frac{r^{n+1}}{(n+1)!}} \otimes \lambda_{n+1}^{-\frac{n(n+1)}{2(n+1)!}r^n} \otimes \lambda_n^{\frac{r^n}{n!}} \otimes \cdots$$
Thus, up to terms of order $r^{2n-1}$,
\begin{eqnarray*}
  \LCH(\X,\LL^r) &=& \lambda_{n+1}^{p(r)r^{n+1}-a_0r^{2n+1}}\otimes \lambda_{n+1}^{ \frac{n(n+1)}{2}a_0r^{2n}}\otimes \lambda_n^{-a_0(n+1)r^{2n}}\otimes \cdots\\
&=&\lambda_{n+1}^{\left(a_1 + \frac{a_0n(n+1)}{2}\right)r^{2n}}\otimes \lambda_n^{-a_0(n+1)r^{2n}}\otimes \cdots\\
&=& \left(\lambda_{n+1}^{\mu + n(n+1)} \lambda_n^{-2(n+1)} \right)^{\frac{a_0}{2}r^{2n}} \otimes \cdots\\
&=& \LCM(\X,\LL)^{\frac{a_0}{2}r^{2n}} \otimes \cdots.
\end{eqnarray*}
where in the second line we have use that $p(r)-a_0r^n = a_1r^{n-1} +
O(r^{n-2})$.  Thus the Proposition follows from this and
\eqref{eq:chexpansion}. Clearly from the above calculation the line
bundle $\epsilon(r,k)$ is of the form claimed.
\end{proof}

\begin{rmk}\label{rmk:modifiedcmline}
  It is convenient to define a $\mathbb Q$-line bundle by
  \begin{equation}
    \LCM'(\X,\LL) = \LCM(\X,\LL)^{\frac{1}{2a_0(n+1)!}},\label{eq:cmmodified}
  \end{equation}
  where $n$ is the relative dimension of $\X$ and $a_0$ is the leading
  order term of the Hilbert polynomial of the fibres.  Then the
  previous Proposition becomes
$$ \LHilb(\X,\LL^r,k) = \LCM'(\X,L)^{a_0^2k^{n+1}r^{2n}} \otimes \cdots.$$ \vspace{4mm}
\end{rmk}
\section{Positivity of the CM line on the stable locus}\label{sec:positivity}

We now make stability assumptions on a general fibre of $\X$ to
deduce positivity properties of $\LCM$.  This idea was first exploited
by Cornalba--Harris to show certain divisor classes on the moduli space of curves are nef 
\cite{cornalba_harris(88):divis_class_assoc_to_famil}.  We recall the
definition of Hilbert stability of a polarised scheme $(X,L)$. Suppose
first that $L$ is very ample.  Then, up to change of coordinates, we
have an embedding of $X$ in $\PP^N$.  For $k\gg
0$ this yields a point $\Hilb(X,L,k)$ in the Hilbert Scheme
$\Hilb(\PP^{N},k)$ representing the subspace of degree $k$ polynomials
on $\PP^{N}$ that vanish along $X$.

\begin{defn}
  We say that a polarised variety $(X,L)$ is \emph{Hilbert semistable}
  if for arbitrarily large $k$ the point $\Hilb(X,L,k)$ is semistable
  in the sense of GIT under the natural $SL_{N+1}$ action.  (This is
  to be understood with respect to the hyperplane bundle on
  $\Hilb(\PP^n,k)$ coming from the embedding in the Grassmannian of
  subspaces of degree $k$ polynomials on $\PP^N$).  We say that
  $(X,L)$ is \emph{asymptotically Hilbert semistable} if $(X,L^r)$ is
  Hilbert semistable for arbitrarily large $r$.
\end{defn}

\begin{thm}[Cornalba--Harris]\label{thm:cornalbaharris}
  Suppose that $\LL$ is relatively very ample and without higher cohomology,
  and that for some $b\in B$ the fibre $(\X_b,\LL_b)$ of $(\X,\LL)$ is
  Hilbert semistable.  Then for arbitrarily large $k$, 
  $\LHilb(\X,\LL,k)$ is effective.
  
  Moreover there is a section of some power of $\LHilb(\X,\LL,k)$ that
  does not vanish at $b$.  Thus if $C\subset B$ is an irreducible
  curve containing $b$ then for arbitrarily large $k$
  $$c_1(\LHilb(\X,\LL,k)).C\geq0.$$
\end{thm}
\begin{proof}[Sketch proof (see \cite{cornalba_harris(88):divis_class_assoc_to_famil} or  \cite{harris-morrison(98):modul_curves}
  Chapter 6D)] For $k\gg 0$ there is a generically surjective map of
  vector bundles
$$\alpha_k\colon S^k \det \pi_*(\LL) \to \pi_*(\LL^k)$$
where $S^k$ denotes the $k$-th symmetric product.  Thus we get an induced map 
$$\beta_k \colon \Lambda^{p(k)}\left(S^k \det \pi_*(\LL)\right) \to \Lambda^{p(k)}\pi_*(\LL^k) = \det \pi_*(\LL^k).$$
The kernel of this map at a point $b\in B$ is the Hilbert point of the
fibre $(\X_b,\LL_b)$.  The stability assumption means that for
arbitrarily large $k$ there is an SL-invariant homogeneous polynomial
$P$ that does not vanish at $\beta_k|_b$.  We may assume that
$P$ has degree $p(1)m$ for some $m\ge 1$.  Thus picking a
trivialisation of $\pi_*(\LL)$ gives a regular local function $f$  on $B$ given
by $f(t) = P(\beta_k|_t)$.  The fact that $P$ is SL-invariant implies
that changing this choice of trivialisation by multiplication by a matrix $M$
scales $f$ by a factor $(\det M)^{-p(k)k/p(1)}$.  Thus $f$ defines a
section of
$$\det \pi_*(\LL^k)^{p(1)m} \otimes \det(\pi_* \LL)^{-kp(k)m}=\LHilb(\X,\LL,k)^m$$
that does not vanish at $b$. It follows that $\LHilb(\X,\LL,k)^m$ has non-negative degree on $C$, finishing the proof.
\end{proof}

\begin{thm}\label{thm:cmlimit}
  Assume that for at least one $b$ the fibre $(\X_b,\LL_b)$ is
  asymptotically Hilbert semistable.  Then $\LCM$ is pseudoeffective.

  Moreover if $C$ is a curve in $B$ and there is a $b\in C$ such that
  $(\X_b,\LL_b)$ is asymptotically Hilbert semistable then
  $c_1(\LCM).C\ge 0$.
\end{thm}
\begin{proof}
  Pick $b$ so $(\X_b,\LL_b)$ is asymptotically Hilbert stable.  This
  means that for arbitrarily large $r$ the fibre $(\X_b,\LL_b^r)$ is
  Hilbert semistable.  By Theorem \ref{thm:cornalbaharris} this
  implies $\LHilb(\X,\LL^r,k)$ is effective for arbitrarily large $k$
  (and has positive degree along any irreducible curve $C$ containing
  $b$).  From \eqref{prop:leadingorder} we know that $\LCM(\X,\LL)$ is
  the leading order term of $\LHilb(\X,\LL^r,k)$.  Thus letting $k$
  tend to infinity and then $r$ tend to infinity we deduce that
  $\LCM(\X,\LL)$ is pseudoeffective (and has non-negative degree along
  $C$).
\end{proof}

%It may be that the CM line is in fact effective in general, and that it
%is free on a suitable ``stable'' locus.

\begin{rmk}\ 
  \begin{enumerate}
  \item In \cite{cornalba_harris(88):divis_class_assoc_to_famil}
    Cornalba--Harris use analogous reasoning to show that
    $\LCH(\X,\LL^r)$ defined in \eqref{eq:linecornalba} is
    pseudoeffective.  Viehweg has a similar argument in
    \cite{viehweg(83):weak_posit_addit_kodair_dimen} to show
    positivity of certain line bundles connected to families of
    general type surfaces.  Ultimately Viehweg relies on Gieseker's
    proof that surfaces of general type are asymptotically Hilbert
    stable with respect to high pluricanonical embeddings
    \cite{gieseker(77):global_modul_for_surfac_gener_type}.
    
  \item That the CM line is the leading order piece of $\LCH$ was
    originally observed by the first author in
    \cite{fine:fibrat_with_const_scalar_curvat}.  The difference
    between this and the argument given there is that we use a
    different definition of the CM line and so avoid any mention of a
    relative canonical bundle of $\X\to B$.
  \end{enumerate}
\end{rmk}

\begin{rmk}\label{cscK rmk}\
	\begin{enumerate}
	\item  Recent work connecting stability to the existence of constant scalar
  curvature K\"ahler metrics gives conditions under which the above theorem can
  be applied.  Suppose there is a $b\in B$ such that the fibre
  $(\X_b,\LL_b)$ is smooth, has no infinitesimal automorphisms, and
  that there exists a K\"ahler metric on $\X_b$ which has constant
  scalar curvature and whose cohomology class is $c_1(\LL_b)$.  Then
  Donaldson shows that $(\X_b,\LL_b)$ is asymptotically Hilbert
  semistable \cite{donaldson(01):scalar_curvat_projec_embed} and thus
  \eqref{thm:cmlimit} applies.  Such a metric exists, for example, when either
    \begin{enumerate}
    \item $\LL_b=K_{\X_b}$ is the canonical bundle (as then a
      K\"ahler--Einstein metric exists by results of Aubin and Yau
      \cite{aub76:equations_monge_ampere,yau78:ricci_kaehler_monge_ampere});
      or
    \item $K_{\X_b}$ is trivial and $\LL$ is arbitrary (by existence of
      a Ricci flat metric by Yau's proof of the Calabi conjecture
      \cite{yau78:ricci_kaehler_monge_ampere}).
    \end{enumerate}
    \item
    A related result is proved by Fujiki--Schumacher \cite{fujiki_schumacher(90):modul_space_extrem_compac_kaehl}. Suppose that $B$ is a smooth curve, $\X \to B$ is a non-trivial family whose fibres are all smooth with no holomorphic vector fields and, moreover, that for each $b$ the ample class $c_1(\LL_b)$ on $\X_b$ contains a constant scalar curvature K\"ahler metric. They prove that $c_1  (\LCM) . B > 0$. Note that this is neither weaker or stronger than the conclusions drawn above. On the one hand, ample implies nef, on the other hand, Theorem \ref{thm:cmlimit} makes no smoothness assumptions. 
    \end{enumerate}
    
\end{rmk}

We can also use Proposition \ref{prop:leadingorder} to prove other facts
about the CM line.  For instance it is immediate that it is rigid, i.e.\ unchanged when $\LL$ is twisted by a line pulled back from the base:
\begin{cor}
 For any line bundle $\sigma$ on $B$,
\begin{eqnarray*}
  \LCM(\X,\LL\otimes\pi^* \sigma) &=& \LCM(\X,\LL) \label{eq:rigid}
%  \LCM'(\X,\LL\otimes\pi^* \sigma) &=& \LCM'(\X,\LL)
%\LCM(\X,\LL^k) &=& \LCM(\X,\LL)^{k^{n+1}} \label{eq:scale}
\end{eqnarray*}
\end{cor}
\begin{proof}
  This follows from \eqref{thm:cmlimit} and the fact that $\LHilb$ is unchanged if $\LL$ is replaced by $\LL\otimes \pi^* \sigma$.
\end{proof}

Another application concerns the CM line on products.  Let $\X_1\to B$ and $\X_2\to B$ be proper flat families of schemes of constant
relative dimension $n$ and $m$ respectively.  For $i=1,2$ let $\LL_i\to
\X_i$ be relatively ample line bundles.  Denote by $\X=\X_1\times_B \X_2\to
B$ the fibred product with projections $q_i\colon \X\to \X_i$, and by
$\LL=q_1^* \LL_1\otimes q_2^*\LL_2$ the product polarisation.

\begin{cor}\label{cor:products}
With notation as above and $\LCM'$ defined as in \eqref{eq:cmmodified},
   $$\LCM'(\X_1\times_B \X_2,\LL)=\LCM'(\X_1,\LL_1) \otimes \LCM'(\X_2,\LL_2).$$ 
 \end{cor}
\begin{proof}
  It is sufficient to prove the formula for the product $\X_1\times
  \X_2\to B\times B$ for then restriction to the diagonal in $B\times
  B$ yields the same formula for the fibred product.  Let $p$ (resp.\ $p_1,p_2$) be the Hilbert polynomial of the fibres of $\X$ (resp.\ $\X_1,\X_2$), so  $p=p_1p_2$. Let $\lambda_\X (k) = \det \pi_! (\LL^k)$ (and similarly for $\lambda_{\X_i}$).   Then by the K\"unneth formula for $k\gg 0$,  $\lambda_{\X}(k) = \lambda_{\X_1}(k)^{p_2(k)}\otimes  \lambda_{\X_2}(k)^{p_1(k)}$ so 
\begin{eqnarray*}
\LHilb(\X,\LL^r,k) &=& \LHilb(\X_1,\LL^r,k)^{p_2(kr) p_2(r)}\otimes \LHilb(\X_2,\LL^r,k)^{p_1(kr)p_1(r)}.
\end{eqnarray*}
Applying \eqref{rmk:modifiedcmline}  and taking the leading order term in $k,r$ gives the result.
\end{proof}

\section{Examples on the unstable locus}

In this section we give examples of families of polarised manifolds
for which the CM line is strictly negative.  Of course in light of
Theorem \ref{thm:cmlimit} each fibre in such a family must be
unstable.  Although it is possible to calculate directly using the definition of the CM line in terms of the $\lambda_i$ it is more transparent to restrict to the case when $\pi\colon \X\to B$ is smooth and $B$ is an irreducible curve.  We will use the first Chern classes of the CM line which, from \eqref{eq:cmchernclass}, is
\begin{equation*}
c_1(\LCM) = \pi_* \left[\mu c_1(\LL)^{n+1} +(n+1)c_1(K_{\X/B})c_1(\LL)^n\right].
\end{equation*}

The examples are obtained by blowing up points in the
fibres of suitable families $\X\to B$.  Suppose $C\subset \X$ is an
irreducible reduced curve dominating $B$ and $\pi\colon C\to B$ has
degree $d$.  Consider the blowup $q\colon \tilde{\X}\to \X$ of $\X$
along $C$ with exceptional divisor $E$.  It is flat over $B$ and generically the fibre $\tilde{\X}_b$ is the fibre $\X_b$
blown up at the $d$ points consisting of the intersection with $C$.

Fix a relatively ample line bundle $\LL\to \X$ and consider the $\mathbb Q$-line
bundle
$$\LL_{\epsilon} = q^*\LL \otimes \OO(-\epsilon E),$$
on $\tilde{\X}$ which is relatively ample for positive $\epsilon$
sufficiently small.  We now  compare the CM lines of $(\X,\LL)$ and $(\tilde{\X},\LL_\epsilon)$.

\begin{prop}\label{prop:blowup}
Suppose $\pi\colon \X\to B$ has relative dimension $n=2$.  Then
  $$c_1(\LCM(\tilde{\X},\LL_\epsilon)) = c_1(\LCM(\X,\LL)) + \epsilon
  \sigma + O(\epsilon^2)$$ where $$\sigma =
  \pi_*\left(-\frac{d}{a_0}c_1(\LL)^3 + 6 c_1(\LL|_C)\right)$$ 
\end{prop}
\begin{proof}
Since $C$
  has codimension $2$, the relative canonical divisors of $\tilde{\X}$
  and $\X$ are related by
$$K_{\tilde{\X}/B} = q^*K_{\X/B} + E.$$
Notice that 
$$c_1(q^*\LL)^2\cdot c_1(E) = q^* c_1(\LL|_C)^2=0$$ as $C$ has dimension 1.  Similarly $c_1(q^*\LL)\cdot c_1(q^*K_{\X/B})\cdot c_1(E)=0$, and as $\OO(E)|_E$ is the tautological line, $-c_1(q^*\LL)\cdot c_1(E)^2 = c_1(\LL|_C)$.  Thus
\begin{eqnarray*}
  c_1(\LL_\epsilon)^3 &=& q^*c_1(\LL)^3 + O(\epsilon^2)\\
  c_1(K_{\tilde \X/B}).c_1(\LL_\epsilon)^2 &=& [c_1(q^*K_{\X/B}) + c_1(E)].[c_1(q^*\LL)^2-2 \epsilon c_1(q^*\LL)c_1(E) + O(\epsilon^2)]\\
  &=& c_1(K_{\X/B}).c_1(\LL)^2 +2\epsilon c_1(\LL|_C) + O(\epsilon^2).
\end{eqnarray*}
Denote the Hilbert polynomial of the fibres of $(\tilde
\X,\LL_\epsilon)$ by $p_\epsilon(k) = a_0(\epsilon)k^2 + a_1(\epsilon)k
+ O(k^0)$, so $p_0(k)=p(k)=a_0k^2+a_1k^1 + O(k^0)$ is the Hilbert polynomial of the fibre of
$(\X,\LL)$.  By the Riemann--Roch theorem applied to
$\tilde{\X}_b$ and $\X_b$,
\begin{eqnarray*}
  a_0(\epsilon) &=& \frac{1}{2} \int_{\tilde{\X}_b} c_1(\LL_\epsilon)^2= \frac{1}{2} \int_{\tilde{\X}_b}(c_1(q^*\LL) - \epsilon c_1(E))^2  =  a_0 + O(\epsilon^2)\\
a_1(\epsilon) &=& -\frac{1}{2}  \int_{\tilde{\X}_b}c_1(K_{\tilde{\X}/B}).c_1(\LL_{\epsilon})=a_1 - \frac{d\epsilon}{2},
%&=&- \frac{1}{2}  \int_{\tilde{\X}_b}(c_1(q^*K_{\X/B}) + c_1(E))\cdot (c_1(\LL) + \epsilon c_1(E)) \\
\end{eqnarray*}
so as $\mu=2a_1/a_0$,
$$\mu_\epsilon := \frac{2a_1(\epsilon)}{a_0(\epsilon)} = \mu -\frac{d\epsilon}{a_0} + O(\epsilon^2).$$
Putting this together, if $\tilde{\pi} = \pi\circ q$ then
\begin{equation}
\begin{split}
c_1&(\LCM(\tilde{\X},\LL_{\epsilon}))=\tilde{\pi}_*[\mu_\epsilon c_1(\LL_\epsilon)^3 + 3c_1(K_{\tilde{\X}/B})c_1(\LL_\epsilon)^2]\\
&= \pi_*[(\mu - \frac{d\epsilon}{a_0})c_1(\LL)^3 + 3(c_1(K_{\X/B})\cdot c_1(\LL)^2 +2\epsilon c_1(\LL)\cdot C) + O(\epsilon^2)]\\
&= c_1(\LCM(\X,\LL)) + \epsilon \sigma  +O(\epsilon^2)
\end{split}
\end{equation}
which proves the Proposition.
\end{proof}

We will apply this to examples where $\LCM(\X,\LL)$ has zero degree and
$\sigma<0$, obtaining examples of families of (smooth) polarised varieties for which the CM line has strictly negative degree.

%The first is similar to
%Morrison's example in
%\cite{cornalba-harris(88):divis_class_assoc_to_famil}.

\begin{example}\label{ex:singleblowup}
  Let $B=\PP^1$ and consider the vector bundle $$\E=\OO_{\PP^1}(2)
  \oplus \OO_{\PP^1}(-1) \oplus \OO_{\PP^1}(-1).$$ We let $\X=\PP(\E)$ (the space of \emph{lines} in $\E$)  and $q\colon \tilde{\X}\to \X$ be the blowup along the curve
  $C=\PP(\OO_{\PP^1}(2))\cong B$.  Set $\LL=\OO_{\PP(\E)}(1)$ which is
  relatively ample and $\LL_\epsilon = q^*\LL \otimes \OO(-\epsilon E)$
  where $E$ is the exceptional set of the blowup.

  As $c_1(\E)=0$ we have for $k\gg 0$, $c_1(\pi_* \LL^k) = c_1(S^k
  \E^*) =0$.  Hence $c_1(\lambda_i)=0$ for all $i$, which implies
  $c_1(\LCM(\X,\LL))=0$.  Moreover, $\pi_*\left(c_1(\LL)^3\right) =0$ and
  $c_1(\LL|_C) = -2$ so $\sigma =-12$.  Thus by Proposition
  \ref{prop:blowup},
  $$c_1(\LCM(\tilde{\X},\LL_\epsilon))<0 \quad \text{for }0<\epsilon \ll 1.$$

  Notice that the fibres of $\tilde{\X}$ consist of $\PP^2$ blown up
  at a single point, which is known to be unstable with respect to any
  polarisation
  (\cite{ross_thomas:obstr_to_exist_const_scalar_curvat_kaehl_metric},
  Example 5.27).  
\end{example}

\begin{rmk}
  This previous example is similar to one due to Morrison
  \cite{cornalba_harris(88):divis_class_assoc_to_famil} showing that
  the Cornalba--Harris line $\LCH$ \eqref{eq:linecornalba} is negative on a family of Steiner
  surfaces.  We remark that in fact $\LCH(\X,\LL^r)$ is negative for
  all $r\gg 0$.
\end{rmk}

\begin{example}\label{ex:positiveunstable}
  In Example \ref{ex:singleblowup} we could instead have let
  $C=\PP(\OO_{\PP^1}(-1))$.  Then the fibres of $\tilde{\X}$ are still
  all unstable with respect to any polarisation.  However $c_1(\LL|_C) =
  1$ and $\sigma = 6$ so $c_1(\LCM(\tilde{\X},\LL_{\epsilon}))>0$.
  Thus the CM line need not always be strictly negative
  on the unstable locus.  Notice that the Hilbert polynomial
  of the fibres in this example are the same as in Example
  \ref{ex:singleblowup}.  Thus there are fixed Hilbert schemes on
  which the CM line is neither strictly positive nor strictly negative
  on the unstable locus.
\end{example}

\begin{example}\label{ex:multipleblowup}
  The fibres of the family in Example \ref{ex:singleblowup} have
  non-trivial infinitesimal automorphisms.  However, this is not
  necessary for the CM line to be negative.  To see this consider the
  family $(\X,\LL_{\epsilon})$ in Example \ref{ex:singleblowup} and fix
  $0<\epsilon\ll 1$ so $\LCM(\tilde{\X},\LL_\epsilon)$ is strictly
  negative.  Pick sections $C_1,C_2,C_3$ of $\pi$ so that
  $C,C_1,C_2,C_3$ meet the generic fibre of $\X$ in 4 generic points.
  Again let $\tilde{\X}$ be the blowup of $\X$ along $C$ and denote
  their proper transforms of $C_1,C_2,C_3$ in $\tilde{\X}$ also by
  $C_i$ (so $C_i$ is disjoint from $E$).

  Let $\Y$ be the blowup of $\tilde{\X}$ along $\cup C_i$, with
  exceptional set $E'$.  For $0<\epsilon'\ll \epsilon$ the line
  $\LL'=\LL \otimes \OO(- \epsilon E)\otimes \OO(- \epsilon 'E')$ is
  relatively ample,  and by continuity $\LCM(\Y,\LL')$ is also strictly
  negative for $0<\epsilon'\ll \epsilon$.

  Notice that a fibre of $\Y$ is $\PP^2$ blown up at 4 distinct
  points, which do not all lie on a line, and thus has discrete
  automorphism group. These manifolds are known to be unstable with
  respect to the polarisations considered here 
  (\cite{ross_thomas:obstr_to_exist_const_scalar_curvat_kaehl_metric},
  Example 5.30)
\end{example}

\begin{example}\label{ex:products}
  Suppose $(\X_1,\LL_1)$ (resp.\ $(\X_2,\LL_2)$) are families over the
  same curve $B$ chosen so that the associated CM line is strictly
  negative (resp.\ non-positive).  Set $\X=\X_1\times_B \X_2$ with
  projections $q_i\colon \X\to \X_1$ and $\LL=q_1^*\LL_1\otimes
  q_2^*\LL_2$.  Then Corollary \ref{cor:products} implies then the CM line
  of $(\X,\LL)$ will have strictly negative degree.  This gives the
  existence of families of manifolds of arbitrary dimension for which
  the CM line is strictly negative.  Notice that since the fibres of
  $\X_1$ are necessarily unstable, the fibres of $\X$ will be as well.
\end{example}

\begin{rmk}
  It is interesting to note that the idea above can, in principle at least, be
  used to show that $\LCM$ is ample.  Still with $n=2$, suppose that
  $\sigma<0$ and
  \begin{itemize}
  \item Some fibre $(\X_b,\LL_b)$ of $(\X,\LL)$ is smooth, has a discrete
    automorphism group and $c_1(\LL_b )$ admits a constant scalar
    curvature K\"ahler metric
  \item The curve $C$ meets $\X_b$ in $d$ distinct points.
  \end{itemize}

  Arrezo-Pacard \cite{arezzo-pacard:blowin_up_desin_const_scalar} show
  that under these assumption the blowup of $X_b$ at the points
  $C.X_b$ admits a cscK metric in classes which make the exceptional
  set sufficiently small.  Thus the fibre of $(\X,\LL_\epsilon)$ over
  $b$ admits a cscK metric.  As mentioned previously, a result of Donaldson
  \cite{donaldson(01):scalar_curvat_projec_embed} then implies that
  this fibre is asymptotically Hilbert semistable and so from
  \eqref{thm:cmlimit}, $c_1(\LCM(\tilde{X},\LL_\epsilon)\ge 0$ for
  $0<\epsilon\ll 1$.  By Proposition \ref{prop:blowup}, the assumption
  that $\sigma<0$ then implies that $c_1(\LCM(\X,\LL))>0$.
  Unfortunately, we do not know of a specific example where this can be
  applied, but wonder nevertheless if this idea could be used to show
  ampleness of the CM line in some cases not covered by Fujiki--Schumacher's result (described in Remark \ref{cscK rmk}).
\end{rmk}

\bibliography{biblio}

{\small \noindent {\tt joel.fine@imperial.ac.uk }} \newline
Department of Mathematics, Imperial College,London SW7 2AZ. UK.\\
{\small \noindent {\tt jaross@math.columbia.edu}} \newline
\noindent Department of Mathematics, Columbia University, New York, NY 10027.
USA. \\

\end{document}